\documentclass[oneside,english]{amsart}
\usepackage{ae}
\usepackage{aecompl}
\usepackage[T1]{fontenc}
\usepackage[latin1]{inputenc}
\usepackage{amssymb}

\makeatletter
 \theoremstyle{plain}
\newtheorem{thm}{Theorem}[section]
  \theoremstyle{plain}
  \newtheorem{lem}[thm]{Lemma}
  \theoremstyle{remark}
  \newtheorem{notation}[thm]{Notation}
  \theoremstyle{plain}
  \newtheorem{cor}[thm]{Corollary}

\usepackage[nativepdf]{hyperref}

\usepackage{ifthen}

\newboolean{Details}

\setboolean{Details}{false}

\makeatother
\begin{document}

\title{R\'{e}nyi Dimension and Gaussian Filtering II}

\author{Terry A. Loring}

\thanks{This work was supported in part by DARPA Contract N00014-03-1-0900.}

\email{loring@math.unm.edu}

\curraddr{Department of Mathematics and Statistics, University of New Mexico,
Albuquerque, NM 87131, USA.}



\urladdr{http://www.math.unm.edu/\textasciitilde{}loring/}

\maketitle
\begin{center}
preprint version
\par\end{center}

\markright {R\'{e}nyi dimension and Gaussian filtering II (preprint version)}\markleft {R\'{e}nyi dimension and Gaussian filtering II (preprint version)} 

\begin{center}
\ifthenelse{\boolean{Details}}{{\Large Fully Detailed Version}}{}
\par\end{center}

\section*{Abstract}

We consider convolving a Gaussian of a varying scale $\epsilon$ against
a Borel measure $\mu$ on Euclidean $\delta$-dimensional space. The
$L^{q}$ norm of the result is differentiable in $\epsilon.$ We calculate
this derivative and show how the upper order of its growth relates
to its lower R\'{e}nyi dimension. We assume $q$ is strictly between
$1$ and $\infty$ and that $\mu$ is finite with compact support.

Consider choosing a sequence $\epsilon_{n}$ of scales for the Gaussians
$g_{\epsilon}(\mathbf{x})=\epsilon^{-\delta}e^{-(|\mathbf{x}|/\epsilon)^{2}}.$
Let $\| f\|_{q}$ denote the $L^{q}$ norm for Lebesgue measure. The
differences\[
\left|\left\Vert g_{\epsilon_{n+1}}\ast\mu\right\Vert _{q}-\left\Vert g_{\epsilon_{n}}\ast\mu\right\Vert _{q}\right|\]
between the norms at adjacent scales $\epsilon_{n}$ and $\epsilon_{n-1}$
can be made to grow more slowly than any positive power of $n$ by
setting the $\epsilon_{n}$ by a power rule. The correct exponent
in the power rule is determined by the lower R\'{e}nyi dimension.

We calculate and find bounds on the derivative of the Gaussian kernel
versions of the correlation integral. We show that a Gaussian Kernel
version of the R\'{e}nyi entropy sum in continuous.

\tableofcontents{}

\section{Differences of Gaussian filters\label{sec:Differences-of-Gaussian}}

Suppose $\mu$ is a finite Borel measure $\mu$ on $\mathbb{R}^{\delta}$
with compact support. If $\delta=2$ we can think of $\mu$ as the
abstraction of an image. In the context of image processing, it is
common to look at the difference of convolutions of two Gaussians
against $\mu.$ This is the case in the standard construction of a
Laplacian pyramid (\cite{BurtAdelson}). Traditionally, the scales
of the kernels vary geometrically. For some purposes, it might be
better to set the scales of the convolution kernels in a different
pattern.

Given a function $g$ on $\mathbb{R}^{\delta},$ thought of as acting
as a canonical filter kernel on measures, we rescale it as \[
g_{\epsilon}(\mathbf{x})=\epsilon^{-\delta}g\left(\epsilon^{-1}\mathbf{x}\right).\]
The most important cases we have in mind are where $g$ is a Gaussian
or a function of compact support that approximates a Gaussian. We
now consider the problem of select some $\epsilon_{n}\searrow0$ so
that the differences\[
g_{\epsilon_{n}}\ast\mu-g_{\epsilon_{n-1}}\ast\mu\]
behave nicely.

We will use the notation $m$ for Lebesgue measure,\begin{eqnarray*}
\| f\|_{q} & = & \left(\int_{\mathbb{R}^{\delta}}f(\mathbf{x})^{q}\, dm(\mathbf{x})\right)^{\frac{1}{q}}\\
 & = & \left(\int_{\mathbb{R}}\int_{\mathbb{R}}\cdots\int_{\mathbb{R}}f(\mathbf{x})^{q}\, dx_{1}dx_{2}\cdots dx_{\delta}\right)^{\frac{1}{q}}.\end{eqnarray*}

One criterion for the selection of the scales $\epsilon_{n}$ is to
keep \[
\left\Vert g_{\epsilon_{n}}\ast\mu-g_{\epsilon_{n-1}}\ast\mu\right\Vert _{q}\]
approximately constant for some choice of $1<q<\infty.$ We don't
see a means of estimating this norm of differences, so we instead
look at the difference of norms.

If $\mu$ is ``fractal,'' then we expect that setting $\epsilon_{n}$
in a geometric series will lead to exponential growth in \[
\left|\left\Vert g_{\epsilon_{n}}\ast\mu\right\Vert _{q}-\left\Vert g_{\epsilon_{n-1}}\ast\mu\right\Vert _{q}\right|.\]
If we set the $\epsilon_{n}$ via a power law, $\epsilon_{n}=n^{-t},$
we can reasonably hope that this difference is more or less constant.

Recall, say from \cite{BarGerTch} or \cite{Olsen}, that the upper
and lower R\'{e}nyi dimensions of $\mu$ for index $q$ are defined
by \[
D_{q}^{\pm}(\mu)=\lim_{\epsilon\rightarrow0}\,_{\textrm{inf }}^{\textrm{sup}}\frac{1}{q-1}\frac{\ln\left(S_{\mu}^{q}(\epsilon)\right)}{\ln(\epsilon)},\]
 where the standard partition function $S_{\mu}^{q}(\epsilon)$ is
taken to be

\begin{equation}
S_{\mu}^{q}(\epsilon)=\sum_{\mathbf{k}\in\mathbb{Z}^{\delta}}\mu(\epsilon\mathbf{k}+\epsilon\mathbb{I})^{q},\label{eq: partition function}\end{equation}
and where $\mathbb{I}$ is a $\delta$-fold product of the unit interval
$[0,1).$

The connection with Gaussian convolution is the formula, due to Gu\'{e}rin,
\[
\limsup_{x\rightarrow\infty}\frac{\ln\left(\left\Vert g_{x^{-1}}\ast\mu\right\Vert _{q}\right)}{\ln(x)}=\frac{q-1}{q}\left(\delta-D_{q}^{-}(\mu)\right).\]
See \cite{Guerin} or \cite[Lemma 2.3]{Loring}. For this formula
to be valid, we need some restriction on $g.$ For simplicity, we
consider the case where $g$ is nonnegative and rapidly decreasing.
By rapidly decreasing, we mean that any order derivative $g^{(\alpha)}$
of $g$ exists and decays at infinity more rapidly than any negative
power of $x.$ Certainly, less is needed. See \cite{Guerin}.

Here is the main result.

\begin{thm}
\label{thm:  growth laws-convolution}Suppose $1<q<\infty.$ Suppose
$g:\mathbb{R}^{\delta}\rightarrow\mathbb{R}$ is nonnegative, nontrivial,
rapidly decreasing and is radially nonincreasing. Suppose $\mu$ is
a finite Borel measure on $\mathbb{R}^{\delta}$ with compact support.
Let $I_{q}(\mu)$ denote the set of positive $t$ for which\[
\forall\alpha>0,\ \lim_{n\rightarrow\infty}\frac{\left|\left\Vert g_{n^{-t}}\ast\mu\right\Vert _{q}-\left\Vert g_{(n-1)^{-t}}\ast\mu\right\Vert _{q}\right|}{n^{\alpha}}=0.\]
If $D_{q}^{-}(\mu)<\delta$ then\[
I_{q}(\mu)=\left(0,\frac{q}{q-1}\left(\delta-D_{q}^{-}(\mu)\right)^{-1}\right],\]
while if $D_{q}^{-}(\mu)=\delta$ then\[
I_{q}(\mu)=(0,\infty).\]

\end{thm}
This is in sharp contrast with what happens if the scales of the kernels
are set to geometric growth.

\begin{thm}
\label{thm: exponential growth laws-convolution}Suppose $q,$ $g$
and $\mu$ are as in Theorem \ref{thm:  growth laws-convolution}.
If $D_{q}^{-}(\mu)<\delta$ then\[
\limsup_{n\rightarrow\infty}\frac{\left|\left\Vert g_{2^{-n}}\ast\mu\right\Vert _{q}-\left\Vert g_{2^{-n+1}}\ast\mu\right\Vert _{q}\right|}{n^{\alpha}}=\infty\quad(\forall\alpha>0).\]

\end{thm}
The proofs require an estimate on the derivative \[
\frac{d}{d\lambda}\ln\left(\left\Vert g_{e^{\lambda}}\ast\mu\right\Vert _{q}\right),\]
which we derive in Section \ref{sec:Differentiating-the-Norm}. Section
\ref{sec:Asymptotic-indices} contains lemmas on the upper order of
positive function of a positive variable and completes the proofs
of the main theorems.

Let us use the notation

\begin{eqnarray*}
\| f\|_{\mu,q} & = & \left(\int_{\mathbb{R}^{\delta}}f(\mathbf{x})^{q}\, d\mu(\mathbf{x})\right)^{\frac{1}{q}}.\end{eqnarray*}
The calculations from Section \ref{sec:Differentiating-the-Norm}
can be adjusted to estimate\[
\frac{d}{d\lambda}\ln\left(\left\Vert g_{e^{\lambda}}\ast\mu\right\Vert _{\mu,q-1}\right).\]
Equivalently, we find bounds on the derivative of\[
\ln\left(\int_{\mathbb{R}^{\delta}}\left(\int_{\mathbb{R}^{\delta}}g\left(\frac{\mathbf{x}-\mathbf{y}}{e^{\lambda}}\right)\, d\mu(\mathbf{y})\right)^{q}\, d\mu(\mathbf{x})\right).\]
This is the content of Section \ref{sec:Gaussian-Kernel-Correlation},
which does not depend on Section \ref{sec:Asymptotic-indices}. This
should be of interest as it relates to computing the correlation dimension
by probabilistic methods. 

In Section \ref{sec:Gaussian-Kernel-Sums} we consider adjusting the
standard partition sum $S_{\mu}^{q}(\epsilon)$ by allowing soft cut-offs
between the cells (bins). We cannot determine the derivative of the
these Gasssian-Kernel sums, but do demonstrate they are continuous
in $\epsilon.$

\section{Differentiating the norm of a filtered measure \label{sec:Differentiating-the-Norm}}

Recall that given a function $g$ on $\mathbb{R}^{\delta}$ we rescale
it as \[
g_{\epsilon}(\mathbf{x})=\epsilon^{-\delta}g\left(\epsilon^{-1}\mathbf{x}\right).\]
This section's goals are to compute\[
\frac{d}{d\epsilon}\left\Vert g_{\epsilon}\ast\mu\right\Vert _{q}\]
and to show that\[
\frac{d}{d\lambda}\ln\left(\left\Vert g_{e^{\lambda}}\ast\mu\right\Vert _{q}\right)\]
is bounded.

\begin{lem}
\label{lem:inner derivative formula}Suppose $g:\mathbb{R}^{\delta}\rightarrow\mathbb{R}$
is differentiable, bounded, and has bounded radial derivative. Let
$h:\mathbb{R}^{\delta}\rightarrow\mathbb{R}$ be the negative of the
radial derivative of $g,$\[
h(\mathbf{x})=-\sum_{j=1}^{\delta}x_{j}\frac{\partial g}{\partial x_{j}}.\]
If $\mu$ is a finite Borel measure on $\mathbb{R}^{\delta}$ then\[
\frac{\partial}{\partial\epsilon}\left[(g_{\epsilon}\ast\mu)(\mathbf{x})\right]=\epsilon^{-1}\left((h_{\epsilon}\ast\mu)(\mathbf{x})-\delta(g_{\epsilon}\ast\mu)(\mathbf{x})\right).\]

\end{lem}
\begin{proof}
For any $\mathbf{w},$\begin{eqnarray*}
\frac{\partial}{\partial\epsilon}g(\epsilon\mathbf{w}) & = & \sum_{j=1}^{\delta}\left(\left.\frac{\partial g}{\partial x_{j}}\right|_{\epsilon\mathbf{w}}\right)\frac{\partial}{\partial\epsilon}\left(\epsilon w_{j}\right)\\
 & = & \sum_{j=1}^{\delta}\left(\left.\frac{\partial g}{\partial x_{j}}\right|_{\epsilon\mathbf{w}}\right)w_{j}\\
 & = & -\epsilon^{-1}h(\epsilon\mathbf{w})\end{eqnarray*}
Suppose $\mathbf{x}$ is fixed. Assume $0<a\leq\epsilon\leq b$ and
$\mathbf{y}\in\mathbb{R}^{\delta}.$ Then

\begin{eqnarray*}
 &  & \frac{\partial}{\partial\epsilon}\left[g_{\epsilon}(\mathbf{x}-\mathbf{y})\right]\\
 &  & =\frac{\partial}{\partial\epsilon}\left[\epsilon^{-\delta}g(\epsilon^{-1}(\mathbf{x}-\mathbf{y})\right]\\
 &  & =\epsilon^{-\delta}(-\epsilon h(\epsilon^{-1}(\mathbf{x}-\mathbf{y}))(-\epsilon^{-2})+(-\delta\epsilon^{-\delta-1})g(\epsilon^{-1}(\mathbf{x}-\mathbf{y}))\\
 &  & =\epsilon^{-1}\left(h_{\epsilon}(\mathbf{x}-\mathbf{y})-\delta g_{\epsilon}(\mathbf{x}-\mathbf{y})\right).\end{eqnarray*}
and so\[
\left|\frac{\partial}{\partial\epsilon}\left[g_{\epsilon}(\mathbf{x}-\mathbf{y})\right]\right|\leq a^{-\delta-1}H+\delta a^{-\delta-1}G,\]
where $G$ and $H$ are bounds on $g$ and $h.$ Since $\mu$ is finite,
$g_{\epsilon}(\mathbf{x}-\mathbf{y})$ is integrable in $\mathbf{y}$.
The Dominated Convergence Theorem gives us \begin{eqnarray*}
\frac{\partial}{\partial\epsilon}\left[(g_{\epsilon}\ast\mu)(\mathbf{x})\right] & = & \frac{\partial}{\partial\epsilon}\int_{\mathbb{R}^{\delta}}g_{\epsilon}(\mathbf{x}-\mathbf{y})\, d\mu(\mathbf{y})\\
 & = & \int_{\mathbb{R}^{\delta}}\epsilon^{-1}\left(h_{\epsilon}(\mathbf{x}-\mathbf{y})-\delta g_{\epsilon}(\mathbf{x}-\mathbf{y})\right)\, d\mu(\mathbf{y})\\
 & = & \epsilon^{-1}\left((h_{\epsilon}\ast\mu)(\mathbf{x})-\delta(g_{\epsilon}\ast\mu)(\mathbf{x})\right).\end{eqnarray*}

\end{proof}
\ifthenelse{\boolean{Details}}{

Pulling the partial derivative past the integral used the standard
result, but with integration done against a general Borel measure
$\mu$. Let us state and prove the result just to be sure.

Suppose $F=F(\mathbf{x},z)$ is defined for $\mathbf{x}\in\mathbb{R}^{\delta}$
and $z\in[a,b].$ Suppose\[
\int_{\mathbb{R}^{\delta}}|F(\mathbf{x},z)|\, d\mu(\mathbf{x})<\infty\]
for each $z.$ Suppose $\frac{\partial F}{\partial z}$ exists. Suppose
there is a nonnegative measurable function $g$ on $\mathbb{R}^{d}$
with\[
\int_{\mathbb{R}^{\delta}}g(\mathbf{x})\, d\mu(\mathbf{x})<\infty\]
and \[
\left|\frac{\partial F}{\partial z}\right|\leq g.\]
The mean value theorem gives us\[
\frac{F(\mathbf{x},z_{1})-F(\mathbf{x},z_{2})}{z_{1}-z_{2}}=\frac{\partial F}{\partial z}(\mathbf{x},\tilde{z})\leq g(\mathbf{x}).\]
By the dominated convergence theorem,\begin{eqnarray*}
\frac{d}{dz}\int_{\mathbb{R}^{\delta}}F(\mathbf{x},z)\, d\mu(\mathbf{x}) & = & \lim_{w\rightarrow z}\int_{\mathbb{R}^{\delta}}\frac{F(\mathbf{x},w)-F(\mathbf{x},z)}{w-z}\, d\mu(\mathbf{x})\\
 & = & \int_{\mathbb{R}^{\delta}}\lim_{w\rightarrow z}\frac{F(\mathbf{x},w)-F(\mathbf{x},z)}{w-z}\, d\mu(\mathbf{x})\\
 & = & \int_{\mathbb{R}^{\delta}}\frac{\partial F}{\partial z}\, d\mu.\end{eqnarray*}
}{}

\begin{notation}
We shall use $x\wedge y$ to denote the minimum of two numbers and
$x\vee y$ for their maximum.
\end{notation}
\begin{thm}
\label{thm:derivative formula}Suppose $1<q<\infty.$ Suppose $g:\mathbb{R}^{\delta}\rightarrow\mathbb{R}$
is rapidly decreasing and let $h:\mathbb{R}^{\delta}\rightarrow\mathbb{R}$
denote the negative of the radial derivative of $g.$ Suppose $g\geq0$
and $h\geq0.$ If $\mu$ is a finite Borel measure on $\mathbb{R}^{\delta}$
with compact support then\[
\frac{d}{d\epsilon}\left\Vert g_{\epsilon}\ast\mu\right\Vert _{q}=\frac{{\displaystyle \int_{\mathbb{R}^{\delta}}\left(g_{\epsilon}\ast\mu\right)^{q-1}(h_{\epsilon}\ast\mu)\, dm}}{\epsilon\left\Vert g_{\epsilon}\ast\mu\right\Vert _{q}^{q-1}}-\frac{\delta\left\Vert g_{\epsilon}\ast\mu\right\Vert _{q}}{\epsilon}\]
and\[
\frac{d}{d\lambda}\ln\left(\left\Vert g_{e^{\lambda}}\ast\mu\right\Vert _{q}\right)=\frac{{\displaystyle \int_{\mathbb{R}^{\delta}}\left(g_{e^{\lambda}}\ast\mu\right)^{q-1}\left(h_{e^{\lambda}}\ast\mu\right)\, dm}}{{\displaystyle \int_{\mathbb{R}^{\delta}}\left(g_{e^{\lambda}}\ast\mu\right)^{q}\, dm}}-\delta.\]

\end{thm}
\begin{proof}
Assume $0<a\leq\epsilon\leq b.$

Pick an integer $k$ with $k>\frac{\delta+1}{q}.$ Since $g$ is rapidly
decreasing there is a $C_{1}$ so that\[
g(\mathbf{x})\leq C_{1}(1\wedge|\mathbf{x}|^{-k}).\]
For all $\mathbf{x},$\begin{eqnarray*}
g_{\epsilon}(\mathbf{x}) & = & \epsilon^{-\delta}g(\epsilon^{-1}\mathbf{x})\\
 & \leq & a^{-\delta}C_{1}(1\wedge(b^{k}|\mathbf{x}|^{-k}))\\
 & = & a^{-\delta}C_{1}b^{k}(b^{-k}\wedge|\mathbf{x}|^{-k}).\end{eqnarray*}
If $|\mathbf{y}|\leq\frac{1}{2}|\mathbf{x}|$ and $2b\leq|\mathbf{x}|$
then\begin{eqnarray*}
g_{\epsilon}(\mathbf{x}-\mathbf{y}) & \leq & a^{-\delta}C_{1}b^{k}(b^{-k}\wedge|\mathbf{x}-\mathbf{y}|^{-k})\\
 & \leq & a^{-\delta}C_{1}b^{k}\left(b^{-k}\wedge\left|\frac{\mathbf{x}}{2}\right|^{-k}\right)\\
 & = & a^{-\delta}C_{1}b^{k}2^{k}|\mathbf{x}|^{-k}.\end{eqnarray*}
Suppose \[
\textrm{supp}(\mu)\subseteq\left\{ \left.\mathbf{y}\in\mathbb{R}^{\delta}\right||\mathbf{y}|\leq R\right\} .\]
If\[
|\mathbf{x}|\geq(2b)\vee(2R)\]
then\begin{eqnarray*}
(g_{\epsilon}\ast\mu)(\mathbf{x}) & = & \int_{|\mathbf{y}|\leq R}g_{\epsilon}(\mathbf{x}-\mathbf{y})\, d\mu(\mathbf{y})\\
 & \leq & \mu\left(\mathbb{R}^{\delta}\right)a^{-\delta}C_{1}b^{k}2^{k}|\mathbf{x}|^{-k}.\end{eqnarray*}
If\[
|\mathbf{x}|\leq(2b)\vee(2R)\]
then we have the estimate\begin{eqnarray*}
(g_{\epsilon}\ast\mu)(\mathbf{x}) & = & \int g_{\epsilon}(\mathbf{x}-\mathbf{y})\, d\mu(\mathbf{y})\\
 & \leq & \int a^{-\delta}C_{1}\, d\mu(\mathbf{y})\\
 & = & \mu\left(\mathbb{R}^{\delta}\right)a^{-\delta}C_{1}.\end{eqnarray*}
For some $C_{2}$ and $C_{3},$\[
(g_{\epsilon}\ast\mu)(\mathbf{x})\leq C_{2}(C_{3}\wedge|\mathbf{x}|^{-k})\]
for all $\mathbf{x}.$

We can repeat the previous argument for $h_{\epsilon}.$ Possibly
increasing $C_{2}$ and $C_{3}$ we can have\[
(h_{\epsilon}\ast\mu)(\mathbf{x})\leq C_{2}(C_{3}\wedge|\mathbf{x}|^{-k}).\]
Therefore\begin{eqnarray*}
\left|\frac{\partial}{\partial\epsilon}\left[(g_{\epsilon}\ast\mu)(\mathbf{x})\right]\right| & = & \left|\epsilon^{-1}\left((h_{\epsilon}\ast\mu)(\mathbf{x})-\delta(g_{\epsilon}\ast\mu)(\mathbf{x})\right)\right|\\
 & \leq & a^{-1}\left((h_{\epsilon}\ast\mu)(\mathbf{x})+\delta(g_{\epsilon}\ast\mu)(\mathbf{x})\right)\\
 & \leq & a^{-1}(\delta+1)C_{2}(C_{3}\wedge|\mathbf{x}|^{-k})\end{eqnarray*}
and\begin{eqnarray*}
\left|\frac{\partial}{\partial\epsilon}\left[\left((g_{\epsilon}\ast\mu)(\mathbf{x})\right)^{q}\right]\right| & = & q\left((g_{\epsilon}\ast\mu)(\mathbf{x})\right)^{q-1}\left|\frac{\partial}{\partial\epsilon}\left[(g_{\epsilon}\ast\mu)(\mathbf{x})\right]\right|\\
 & \leq & qa^{-1}(\delta+1)(C_{2}(C_{3}\wedge|\mathbf{x}|^{-k}))^{q}\\
 & = & qa^{-1}(\delta+1)C_{2}^{q}(C_{3}^{q}\wedge|\mathbf{x}|^{-qk}),\end{eqnarray*}
which is is integrable since $k$ is larger than $\frac{\delta+1}{q}.$
We can use dominated convergence again. By Lemma \ref{lem:inner derivative formula},\begin{eqnarray*}
 &  & \frac{d}{d\epsilon}\left(\left\Vert g_{\epsilon}\ast\mu\right\Vert _{q}^{q}\right)\\
 &  & =\int_{\mathbb{R}^{\delta}}q\left((g_{\epsilon}\ast\mu)(\mathbf{x})\right)^{q-1}\epsilon^{-1}\left((h_{\epsilon}\ast\mu)(\mathbf{x})-\delta(g_{\epsilon}\ast\mu)(\mathbf{x})\right)\, dm(\mathbf{x})\\
 &  & =\frac{q}{\epsilon}\left(\int_{\mathbb{R}^{\delta}}\left((g_{\epsilon}\ast\mu)(\mathbf{x})\right)^{q-1}(h_{\epsilon}\ast\mu)(\mathbf{x})\, dm(\mathbf{x})-\delta\int_{\mathbb{R}^{\delta}}\left((g_{\epsilon}\ast\mu)(\mathbf{x})\right)^{q}\, dm(\mathbf{x})\right).\end{eqnarray*}

The two derivatives formulas in the statement of the lemma now follow
from \[
\frac{d}{d\epsilon}\left\Vert g_{\epsilon}\ast\mu\right\Vert _{q}=\frac{1}{q}\left\Vert g_{\epsilon}\ast\mu\right\Vert _{p}^{1-q}\frac{d}{d\epsilon}\left(\left\Vert g_{\epsilon}\ast\mu\right\Vert _{q}^{q}\right)\]
and\[
\frac{d}{d\lambda}\ln\left(\left\Vert g_{e^{\lambda}}\ast\mu\right\Vert _{q}\right)=\frac{\left(\left.\frac{d}{d\epsilon}\right|_{\epsilon=e^{\lambda}}\left\Vert g_{\epsilon}\ast\mu\right\Vert _{q}\right)e^{\lambda}}{\left\Vert g_{e^{\lambda}}\ast\mu\right\Vert _{q}}.\]

\end{proof}
\begin{thm}
Suppose $1<q<\infty.$ Suppose $g:\mathbb{R}^{\delta}\rightarrow\mathbb{R}$
is rapidly decreasing and let $h:\mathbb{R}^{\delta}\rightarrow\mathbb{R}$
denote the negative of the radial derivative of $g.$ Suppose $g\geq0$
and $h\geq0.$ If $\mu$ is a finite Borel measure on $\mathbb{R}^{\delta}$
with compact support then\[
-\delta\epsilon^{-1}\left\Vert g_{\epsilon}\ast\mu\right\Vert _{q}\leq\frac{d}{d\epsilon}\left\Vert g_{\epsilon}\ast\mu\right\Vert _{q}\leq\epsilon^{-1}\| h_{\epsilon}\ast\mu\|_{q}-\delta\epsilon^{-1}\left\Vert g_{\epsilon}\ast\mu\right\Vert _{q}\]
and\[
-\delta\leq\frac{d}{d\lambda}\ln\left(\left\Vert g_{e^{\lambda}}\ast\mu\right\Vert _{q}\right)=\frac{\| h_{e^{\lambda}}\ast\mu\|_{q}}{\| g_{e^{\lambda}}\ast\mu\|_{q}}-\delta.\]

\end{thm}
\begin{proof}
The two lower bounds follow trivially from the last lemma and the
fact that $g$ and $h$ are nonnegative.

H\"{o}lder's inequality gives us\begin{eqnarray*}
 &  & \int_{\mathbb{R}^{\delta}}\left((g_{\epsilon}\ast\mu)(\mathbf{x})\right)^{q-1}(h_{\epsilon}\ast\mu)(\mathbf{x})\, dm(\mathbf{x})\\
 &  & \leq\left(\int_{\mathbb{R}^{\delta}}\left((g_{\epsilon}\ast\mu)(\mathbf{x})\right)^{q}\, dm(\mathbf{x})\right)^{\frac{q-1}{q}}\left(\int_{\mathbb{R}^{\delta}}\left((h_{\epsilon}\ast\mu)(\mathbf{x})\right)^{q}\, dm(\mathbf{x})\right)^{\frac{1}{q}}\\
 &  & =\| g_{\epsilon}\ast\mu\|_{q}^{q-1}\| h_{\epsilon}\ast\mu\|_{q}\end{eqnarray*}
and the upper bounds follows.
\end{proof}
\begin{cor}
\label{cor:  bounded slope of log-log}Suppose $1<q<\infty.$ Suppose
$g:\mathbb{R}^{\delta}\rightarrow\mathbb{R}$ is nonnegative, nontrivial,
rapidly decreasing and is radially nonincreasing. There is a finite
constant $C$ so that if $\mu$ is a finite Borel measure on $\mathbb{R}^{\delta}$
with compact support then\[
-\delta\leq\frac{d}{d\lambda}\ln\left(\left\Vert g_{e^{\lambda}}\ast\mu\right\Vert _{q}\right)\leq C\]
for all $\lambda.$ If $g$ is a Gaussian, then we may take $C=0.$
\end{cor}
\begin{proof}
We can apply \cite[Lemma 2.1]{Loring} to $g.$ With $S_{\mu}^{q}(\epsilon)$
is the partition function (\ref{eq: partition function}), this tells
us that there is a $D$ so that\[
D^{-1}\leq\frac{\epsilon^{\delta\frac{q-1}{q}}\| g_{\epsilon}\ast\mu\|_{q}}{\left(S_{\mu}^{q}(\epsilon)\right)^{\frac{1}{q}}}\leq D.\]
From the proof of \cite[Lemma 2.1]{Loring} we see that $D$ can be
taken to depend only on $q$ and $g.$ This conclusion is valid for
the negative radial derivative $h$ as well. This is because $\| h_{\epsilon}\ast\mu\|_{q}$
is invariant under translations of $h$ and clearly $h$ is bounded
away from zero on some open set. Therefore\[
\frac{\| h_{\epsilon}\ast\mu\|_{q}}{\| g_{\epsilon}\ast\mu\|_{q}}=\frac{\epsilon^{\delta\frac{q-1}{q}}\| h_{\epsilon}\ast\mu\|_{q}}{\left(S_{\mu}^{q}(\epsilon)\right)^{\frac{1}{q}}}\frac{\left(S_{\mu}^{q}(\epsilon)\right)^{\frac{1}{q}}}{\epsilon^{\delta\frac{q-1}{q}}\| g_{\epsilon}\ast\mu\|_{q}}\]
is bounded above and away from zero.

In the Gaussian case, we know that $\| g_{\epsilon}\ast\mu\|_{q}$
is non-increasing (c.f.~\cite[Lemma 3.1]{Loring}) and so the derivative
is nonpositive.
\end{proof}

\section{Asymptotic indices\label{sec:Asymptotic-indices}}

Given a function $f>0$ on the positive reals, the quantities\[
\overline{d}(f)=\limsup_{x\rightarrow\infty}\frac{\ln\left(f(x)\right)}{\ln(x)}\]
and\[
\underline{d}(f)=\liminf_{x\rightarrow\infty}\frac{\ln\left(f(x)\right)}{\ln(x)}\]
are the upper and lower orders of $f.$ These provide a simple way
to compare the asymptotic behavior of $f(x)$ to $x^{c}$ for various
powers $c.$ Equivalently, $\overline{d}(f)$ is the smallest extended
real number so that\begin{equation}
c>\overline{d}(f)\implies f(x)\leq x^{c}\textrm{ for large }x\label{eq: upper order via bounds}\end{equation}
and $\underline{d}(f)$ is the largest extended real number so that\begin{equation}
c<\underline{d}(f)\implies f(x)\geq x^{c}\textrm{ for large }x.\label{eq: lower order via bounds}\end{equation}
The proof is not complicated. See (\cite{Hoscneider}). For a look
at how upper and lower order relate to regular variation, see \cite{BingGoldTeu}.

In broad terms, if\[
c=\lim_{x\rightarrow\infty}\frac{\ln\left(f(x)\right)}{\ln(x)}\]
exists, then $f(x)$ behaves not so differently from $x^{c}.$ There
is no reason to think that if $f^{\prime}(x)$ exists it must behave
like $x^{c-1}.$ However, if there are some bounds on the derivative
of the log-log plot of $f$ then we are able to deduce the upper order
of $|f^{\prime}|$ from the upper order of $f.$

For the lower order on $|f^{\prime}|,$ we have found no particularly
interesting result that can be applied to $\| g_{x^{-1}}\ast\mu\|_{q}.$
The difficulty is that even if we assume $g$ is a Gaussian we don't
know if the derivative of $\| g_{x^{-1}}\ast\mu\|_{q}$ is bounded
away from zero.

We take the liberty of setting $\ln(0)=-\infty,$ and indeed $\frac{\ln(0)}{C}=-\infty.$
(This is to accommodate $f^{\prime}(x)=0$ at some $x$ and $f(n)=f(n-1)$
at some $n.$) Both (\ref{eq: upper order via bounds}) and (\ref{eq: lower order via bounds})
remain valid.

\begin{lem}
Suppose\label{lemma:  discrete limsup ok given bound}\[
f:[1,\infty)\rightarrow(0,\infty)\]
 is differentiable and that for some finite constant $C,$\begin{equation}
\left|\frac{d}{dx}\ln\left(f\left(e^{x}\right)\right)\right|\leq C.\label{eq: bound on derivative of log-log plot}\end{equation}
 Given any nondecreasing sequence $x_{n}$ with limit $\infty,$ if\begin{equation}
\lim_{n\rightarrow\infty}\frac{\ln(x_{n+1})}{\ln(x_{n})}=1\label{eq:slow enough sequence}\end{equation}
then\[
\limsup_{n\rightarrow\infty}\frac{\ln\left(f(x_{n})\right)}{\ln(x_{n})}=\limsup_{x\rightarrow\infty}\frac{\ln\left(f(x)\right)}{\ln(x)}.\]

\end{lem}
\begin{proof}
The bound (\ref{eq: bound on derivative of log-log plot}) implies\[
\left|\ln\left(f\left(e^{x}\right)\right)-\ln\left(f\left(e^{y}\right)\right)\right|\leq C\left|x-y\right|\]
or\[
\left|\ln\left(f\left(x\right)\right)-\ln\left(f\left(y\right)\right)\right|\leq C\left|\ln(x)-\ln(y)\right|.\]
\ifthenelse{\boolean{Details}}{Also \[
\left|\ln\left(f\left(x\right)\right)\right|\leq C\left|\ln(x)\right|+\left|\ln\left(f\left(1\right)\right)\right|.\]

Suppose $x_{n}$ is an nondecreasing sequence and that (\ref{eq:slow enough sequence})
is true. The inequality\[
\limsup_{n\rightarrow\infty}\frac{\ln\left(f(x_{n})\right)}{\ln(x_{n})}\leq\limsup_{x\rightarrow\infty}\frac{\ln\left(f(x)\right)}{\ln(x)}\]
follows from the usual business about subnets.

Choose a sequence $y_{n}\rightarrow\infty$ with\[
\lim_{n\rightarrow\infty}\frac{\ln\left(f(y_{n})\right)}{\ln(y_{n})}=\limsup_{x\rightarrow\infty}\frac{\ln\left(f(x)\right)}{\ln(x)}.\]
Without loss of generality, $x_{n}\geq e$ and $y_{n}\geq e$ for
all $n.$ For each $n$ there is some $m_{n}$ with\[
x_{m_{n}}\leq y_{n}\leq x_{m_{n}+1}.\]
If we set \[
z_{n}=x_{m_{n}}\]
then\[
\lim_{n\rightarrow\infty}\frac{\ln(y_{n})}{\ln(z_{n})}=1.\]
and\[
\limsup_{n\rightarrow\infty}\frac{\ln\left(f(z_{n})\right)}{\ln(z_{n})}\leq\limsup_{n\rightarrow\infty}\frac{\ln\left(f(x_{n})\right)}{\ln(x_{n})}.\]
Finally,\begin{eqnarray*}
 &  & \lim_{n\rightarrow\infty}\left|\frac{\ln\left(f(z_{n})\right)}{\ln(z_{n})}-\frac{\ln\left(f(y_{n})\right)}{\ln(y_{n})}\right|\\
 &  & \leq\lim_{n\rightarrow\infty}\left|\frac{\ln\left(f(z_{n})\right)}{\ln(z_{n})}\right|\left|\frac{\ln(y_{n})-\ln(z_{n})}{\ln(y_{n})}\right|+\frac{\left|\ln\left(f(z_{n})\right)-\ln\left(f(y_{n})\right)\right|}{\ln(y_{n})}\\
 &  & \leq\lim_{n\rightarrow\infty}\frac{C\left|\ln(z_{n})\right|+\left|\ln\left(f\left(1\right)\right)\right|}{\ln(z_{n})}\left|1-\frac{\ln(z_{n})}{\ln(y_{n})}\right|+\frac{C\left|\ln(z_{n})-\ln(y_{n})\right|}{\ln(y_{n})}\\
 &  & \leq\lim_{n\rightarrow\infty}\left(C+\left|\ln\left(f\left(1\right)\right)\right|\right)\left|1-\frac{\ln(z_{n})}{\ln(y_{n})}\right|+C\left|\frac{\ln(z_{n})}{\ln(y_{n})}-1\right|\\
 &  & =0\end{eqnarray*}
so that\begin{eqnarray*}
\limsup_{n\rightarrow\infty}\frac{\ln\left(f(z_{n})\right)}{\ln(z_{n})} & = & \limsup_{n\rightarrow\infty}\frac{\ln\left(f(y_{n})\right)}{\ln(y_{n})}\\
 & = & \limsup_{x\rightarrow\infty}\frac{\ln\left(f(x)\right)}{\ln(x)}.\end{eqnarray*}
}{The rest of the proof mimics that of \cite[Lemma 4.1]{Loring},
and is omitted.}
\end{proof}
\begin{lem}
\label{lemma:  power law of derivative}Suppose\[
f:[1,\infty)\rightarrow(0,\infty)\]
is differentiable. If, for some finite constant $C,$\[
\left|\frac{d}{dx}\ln\left(f\left(e^{x}\right)\right)\right|\leq C\]
for all $x,$ then\begin{eqnarray*}
\limsup_{n\rightarrow\infty}\frac{\ln\left|f(n)-f(n-1)\right|}{\ln(n)} & = & \limsup_{x\rightarrow\infty}\frac{\ln\left|f^{\prime}(x)\right|}{\ln(x)}\\
 & = & \limsup_{x\rightarrow\infty}\frac{\ln\left(f(x)\right)}{\ln(x)}-1.\end{eqnarray*}

\end{lem}
\begin{proof}
Suppose $f$ is a function with the bounds $\pm C$ on the slope of
its log-log plot. For each $n,$ the Mean Value Theorem gives us a
number $x_{n}$ in the range\[
n-1\leq x_{n}\leq n\]
for which\[
f(n)-f(n-1)=f^{\prime}(x_{n}).\]

From basic facts about nets we obtain \begin{eqnarray*}
\limsup_{n\rightarrow\infty}\frac{\ln\left|f(n)-f(n-1)\right|}{\ln(n)} & = & \limsup_{n\rightarrow\infty}\frac{\ln\left|f^{\prime}(x_{n})\right|}{\ln(n)}\\
 & \leq & \limsup_{n\rightarrow\infty}\frac{\ln\left|f^{\prime}(x_{n})\right|}{\ln(x_{n})}\\
 & \leq & \limsup_{x\rightarrow\infty}\frac{\ln\left|f^{\prime}(x)\right|}{\ln(x)}.\end{eqnarray*}
Let \[
g(x)=\ln\left(f\left(e^{x}\right)\right)\]
so that\[
|g^{\prime}(x)|\leq C\]
and\[
g^{\prime}(x)=\frac{e^{x}f^{\prime}\left(e^{x}\right)}{f\left(e^{x}\right)}.\]
This can be rewritten as\[
f^{\prime}(x)=g^{\prime}(\ln(x))x^{-1}f(x)\]
and so we have\[
\left|f^{\prime}(x)\right|\leq Cx^{-1}f(x).\]
Therefore\begin{eqnarray*}
\limsup_{x\rightarrow\infty}\frac{\ln\left|f^{\prime}(x)\right|}{\ln(x)} & \leq & \limsup_{x\rightarrow\infty}\frac{\ln C-\ln(x)+\ln(f(x))}{\ln(x)}\\
 & = & \limsup_{x\rightarrow\infty}\frac{\ln(f(x))}{\ln(x)}-1.\end{eqnarray*}

To finish, we must show\[
\limsup_{x\rightarrow\infty}\frac{\ln(f(x))}{\ln(x)}-1\leq\limsup_{n\rightarrow\infty}\frac{\ln\left|f(n)-f(n-1)\right|}{\ln(n)}.\]
We can apply Lemma \ref{lemma:  discrete limsup ok given bound},
because\[
1\leq\frac{\ln(x_{n+1})}{\ln(x_{n})}\leq\frac{\ln(n+1)}{\ln(n-1)}\rightarrow1,\]
and this tells us that it will suffice to show\[
\limsup_{n\rightarrow\infty}\frac{\ln\left(f(n)\right)}{\ln(n)}-1\leq\limsup_{n\rightarrow\infty}\frac{\ln\left|f(n)-f(n-1)\right|}{\ln(n)}.\]

Let\[
m=\limsup_{n\rightarrow\infty}\frac{\ln\left|f(n)-f(n-1)\right|}{\ln(n)}.\]
Suppose we are given $\delta>0.$ Then pick $c\neq-1$ with\[
m<c<m+\delta.\]
(If $m=\infty$ we have nothing to prove. If $m=-\infty$ then modify
this to picking $c\neq-1$ less then any given finite number $C.$)
There is a natural number $n_{0}$ so that\[
n\geq n_{0}\implies|f(n)-f(n-1)|\leq n^{c}.\]
 For large $n,$\begin{eqnarray*}
f(n) & = & f(n_{0})+\sum_{k=n_{0}+1}^{n}\left|f(k)-f(k-1)\right|\\
 & \leq & f(n_{0})+\sum_{k=n_{0}+1}^{n}k^{c}\\
 & \leq & f(n_{0})+\int_{n_{0}}^{n+1}y^{c}\, dy\\
 & \leq & f(n_{0})+\frac{1}{c+1}(n+1)^{c+1}.\end{eqnarray*}
For large $n,$\[
f(n)\leq n^{m+\delta+1}.\]
Therefore\[
\limsup_{n\rightarrow\infty}\frac{\ln\left(f(n)\right)}{\ln(n)}\leq m+\delta+1.\]
Since this is true for all $\delta>0,$ we are done. (If $m=-\infty$
then we obtain this $\limsup$ is less than $C,$ for all finite $C,$
and so is also $-\infty.$)
\end{proof}
\begin{lem}
\label{lem:power law of adjusted difference}Suppose\[
f:[1,\infty)\rightarrow(0,\infty)\]
is differentiable and that there is a finite constant $C$ so that\[
\left|\frac{d}{dx}\ln\left(f\left(e^{x}\right)\right)\right|\leq C\]
for all $x.$ If \[
\overline{d}(f)=\limsup_{x\rightarrow\infty}\frac{\ln(f(x))}{\ln(x)}>0\]
then\[
\left\{ t>0\left|\forall\alpha>0,\lim_{n\rightarrow\infty}\frac{\left|f(n^{t})-f((n-1)^{t})\right|}{n^{\alpha}}=0\right.\right\} =\left(0,\overline{d}(f)^{-1}\right].\]
If $\overline{d}(f)=0$ then \[
\left\{ t>0\left|\forall\alpha>0,\lim_{n\rightarrow\infty}\frac{\left|f(n^{t})-f((n-1)^{t})\right|}{n^{\alpha}}=0\right.\right\} =(0,\infty).\]

\end{lem}
\begin{proof}
For a sequence $a_{n}>0,$ it is routine to show that\begin{equation}
\limsup_{n}\frac{\ln(a_{n})}{\ln(n)}\leq0\iff\forall\alpha>0,\lim_{n\rightarrow\infty}\frac{a_{n}}{n^{\alpha}}=0.\label{eq:upper order less zero}\end{equation}

Let $h(x)=\ln(f(e^{x}))$ so that \[
|h^{\prime}(x)|\leq C.\]
With $t>0$ to be specified below, define\[
g(x)=f(x^{t}).\]
Then\[
\left|\frac{d}{dx}\ln\left(g\left(e^{x}\right)\right)\right|\leq tC\]
and we may apply Lemma \ref{lemma:  power law of derivative} to $g.$

As to the upper order:\begin{eqnarray*}
\limsup_{x\rightarrow\infty}\frac{\ln(g(x))}{\ln(x)} & = & \limsup_{x\rightarrow\infty}\frac{\ln(f(x^{t}))}{\ln(x)}\\
 & = & \limsup_{x\rightarrow\infty}\frac{\ln(f(x))}{\ln(x^{\frac{1}{t}})}\\
 & = & t\overline{d}(f)\end{eqnarray*}
By Lemma \ref{lemma:  power law of derivative},\begin{eqnarray*}
\limsup_{n\rightarrow\infty}\frac{\ln\left(\left|f(n^{t})-f((n-1)^{t})\right|\right)}{\ln(n)} & = & \limsup_{n\rightarrow\infty}\frac{\ln\left(\left|g(n)-g(n-1)\right|\right)}{\ln(n)}\\
 & = & \limsup_{x\rightarrow\infty}\frac{\ln\left(g(x)\right)}{\ln(x)}-1\\
 & = & t\overline{d}(f)-1\end{eqnarray*}

If $t>\overline{d}(f)$ then by (\ref{eq:upper order less zero})
there exists $\alpha>0$ so that\[
\frac{f(n^{t})-f((n-1)^{t})}{n^{\alpha}}\not\rightarrow0.\]
If $t\leq\overline{d}(f)$ then\[
\frac{f(n^{t})-f((n-1)^{t})}{n^{\alpha}}\rightarrow0\]
for all positive $\alpha.$
\end{proof}
Theorem \ref{thm:  growth laws-convolution} now follows, since if\[
f(x)=\left\Vert g_{x^{-1}}\ast\mu\right\Vert _{q}\]
then by \cite{Guerin}, or \cite[Lemma 2.3]{Loring},\[
\overline{d}(f)=\frac{q-1}{q}\left(\delta-D_{q}^{-}(\mu)\right).\]

To prove Theorem \ref{thm: exponential growth laws-convolution} requires
only the following lemma. 

\begin{lem}
Suppose\[
f:[1,\infty)\rightarrow(0,\infty)\]
is differentiable. If, for some finite constant $C,$\[
\left|\frac{d}{dx}\ln\left(f\left(e^{x}\right)\right)\right|\leq C\]
for all $x,$ and if\[
\overline{d}(f)=\limsup_{x\rightarrow\infty}\frac{\ln(f(x))}{\ln(x)}>0,\]
then\[
\limsup_{n\rightarrow\infty}\frac{\ln\left|f(2^{n})-f(2^{n-1})\right|}{\ln(n)}=\infty.\]

\end{lem}
\begin{proof}
Suppose for some $\alpha>0$ there is an $n_{0}$ so that \[
n\geq n_{0}\implies\left|f(2^{n})-f(2^{n-1})\right|\leq n^{\alpha}.\]
Then \[
n\geq n_{0}\implies f(2^{n})\leq f\left(2^{n_{0}}\right)+n^{\alpha+1}.\]
Suppose $\beta>0.$ Then for some $n_{1}\geq n_{0},$\[
n\geq n_{1}\implies f\left(2^{n_{0}}\right)+n^{\alpha+1}\leq2^{n\beta}.\]
Therefore\[
\limsup_{n\rightarrow\infty}\frac{\ln\left(f\left(2^{n}\right)\right)}{\ln\left(2^{n}\right)}\leq\beta.\]
Lemma \ref{lemma:  discrete limsup ok given bound} tells us\[
\limsup_{n\rightarrow\infty}\frac{\ln\left(f\left(x\right)\right)}{\ln\left(x\right)}=\limsup_{n\rightarrow\infty}\frac{\ln\left(f\left(2^{n}\right)\right)}{\ln\left(2^{n}\right)}=0.\]

\end{proof}

\section{Gaussian kernel correlation integrals\label{sec:Gaussian-Kernel-Correlation}}

The probabilistic interpretations of the correlation integral\[
\int\mu(\mathbf{x}+\epsilon\mathbb{B})^{q-1}\, d\mu(\mathbf{x})\]
make it a common tool for determining the R\'{e}nyi dimensions of
$\mu.$ Here\[
\mathbb{B}=\left\{ \left.\mathbf{x}\in\mathbb{R}^{\delta}\right||\mathbf{x}|\leq1\right\} \]
and $\mu$ is a Borel probability measure on $\mathbb{R}^{\delta}.$
When $q=2$ the correlation integral is \[
\int\mu(B_{\epsilon}(\mathbf{x}))\, d\mu(\mathbf{x})=\mbox{Pr}\left\{ \left|X_{1}-X_{2}\right|\leq\epsilon\right\} ,\]
where $X_{1}$ and $X_{2}$ are random locations in the probability
space $\left(\mathbb{R}^{\delta},\mu\right).$ Using a sharp cut-off
for the allowed distance seems unwise in a numerical situation, as
is discussed in \cite{DiksEstNoisyAttr,DiksCorrelDimReturnsVolatil,GhezVaientiI,GhezVaientiII,Guerin,AlexPangBount,NolteGuidoZieheKlaus,SchreiberInflGaussNoise,Schreiber}.

Consider the expectation of the scalar-valued random variable\[
G\left(\epsilon^{-1}\left|X_{1}-X_{2}\right|\right)\]
for a function such as a Gaussian $G(x)=e^{-x^{2}},$ or any $G\geq0$
that is positive at $0$ and rapidly decreasing. This expectation
can be rewritten as follows. Let $g(\mathbf{x})=G(|\mathbf{x}|)$
and \[
g_{\epsilon}(\mathbf{x})=\epsilon^{-\delta}G(\epsilon^{-1}|\mathbf{x}|).\]
 Then\begin{eqnarray*}
E\left[G\left(\epsilon^{-1}\left\Vert X_{1}-X_{2}\right\Vert \right)\right] & = & \epsilon^{\delta}\int\int g_{\epsilon}\left(\mathbf{x}-\mathbf{y}\right)\, d\mu(\mathbf{y})\, d\mu(\mathbf{x})\\
 & = & \epsilon^{\delta}\left\Vert g_{\epsilon}\ast\mu\right\Vert _{\mu,1}.\end{eqnarray*}

As is shown in \cite{BarGerTch}, \begin{eqnarray*}
D_{q}^{\pm}(\mu) & = & \lim_{\epsilon\rightarrow0}\,_{\textrm{inf }}^{\textrm{sup}}\frac{1}{q-1}\frac{\ln\left(\int_{\mathbb{R}^{\delta}}\left(\epsilon^{\delta}g_{\epsilon}\ast\mu\right)^{q-1}\, d\mu\right)}{\ln(\epsilon)}\\
 & = & \delta+\lim_{\epsilon\rightarrow0}\,_{\textrm{inf }}^{\textrm{sup}}\frac{\ln\left(\left\Vert g_{\epsilon}\ast\mu\right\Vert _{\mu,q-1}\right)}{\ln(\epsilon)}\end{eqnarray*}
and more specifically there is a constant $C\neq0$ so that\[
C^{-1}\leq\frac{\int_{\mathbb{R}^{\delta}}\epsilon^{\delta}\left(g_{\epsilon}\ast\mu\right)^{q-1}\, d\mu}{S_{\mu}^{p}(\epsilon)}\leq C\]
for all $\epsilon.$

The R\'{e}nyi dimensions of $\mu$ can be computed as\[
\lim_{\lambda\rightarrow-\infty}\,_{\textrm{inf }}^{\textrm{sup}}\frac{\ln(P_{\mu}(e^{\lambda}))}{\lambda}\]
for at least the following six choices of partition function. (See
\cite{BarGerTch,Guerin,GuysinskyYaskolko,Pesin} and Section \ref{sec:Gaussian-Kernel-Sums}.)
\begin{equation}
P_{\mu}^{q}(\epsilon)=\left(S_{\mu}^{q}(\epsilon)\right)^{\frac{1}{q-1}}=\left(\sum_{\mathbf{j}\in\mathbb{Z}^{\delta}}\mu(\epsilon\mathbf{j}+\epsilon\mathbb{I})^{q}\right)^{\frac{1}{q-1}};\label{eq:partition_one}\end{equation}
\begin{equation}
P_{\mu}^{q}(\epsilon)=\left(\int_{\mathbb{R}^{\delta}}\mu(\mathbf{x}+\epsilon\mathbb{B})^{q-1}\, d\mu(\mathbf{x})\right)^{\frac{1}{q-1}};\label{eq:partition_two}\end{equation}
\begin{equation}
P_{\mu}^{q}(\epsilon)=\left(\int_{\mathbb{R}^{\delta}}\mu(\mathbf{x}+\epsilon\mathbb{B})^{q}\frac{1}{\epsilon^{\delta}}\, dm(\mathbf{x})\right)^{\frac{1}{q-1}};\label{eq:partion_three}\end{equation}
\begin{equation}
P_{\mu}^{q}(\epsilon)=\left(\sum_{\mathbf{j}\in\mathbb{Z}^{\delta}}\left(\int_{\mathbb{R}^{\delta}}g\left(\mathbf{j}-\frac{\mathbf{y}}{\epsilon}\right)\, d\mu(\mathbf{y})\right)^{q-1}\right)^{\frac{1}{q-1}};\label{eq:partition_four}\end{equation}
\begin{equation}
P_{\mu}^{q}(\epsilon)=\left(\int_{\mathbb{R}^{\delta}}\left(\int_{\mathbb{R}^{\delta}}g\left(\frac{\mathbf{x}-\mathbf{y}}{\epsilon}\right)\, d\mu(\mathbf{y})\right)^{q-1}\, d\mu(\mathbf{x})\right)^{\frac{1}{q-1}};\label{eq:partition_five}\end{equation}
\begin{equation}
P_{\mu}^{q}(\epsilon)=\left(\int_{\mathbb{R}^{\delta}}\left(\int_{\mathbb{R}^{\delta}}g\left(\frac{\mathbf{x}-\mathbf{y}}{\epsilon}\right)\, d\mu(\mathbf{y})\right)^{q}\frac{1}{\epsilon^{\delta}}\, dm(\mathbf{x})\right)^{\frac{1}{q-1}}.\label{eq:Partion_six}\end{equation}
The functions (\ref{eq:partition_one}) and (\ref{eq:partition_two})
can be discontinuous in $\epsilon.$ (A sum of two point masses shows
this.) Assuming $\mu$ has compact support, we find that (\ref{eq:partition_four})
is continuous for $1<q<\infty$ (Theorem \ref{thm:continuityGaussKerSums}),
that (\ref{eq:partition_five}) continuous for $1<q<\infty$ and differentiable
for $2<q<\infty,$ (Theorems \ref{thm:continuityCorrelationInt} and
\ref{thm:derivative formula 2}), and that (\ref{eq:Partion_six})
is differentiable for $1<q<\infty$ (Theorem \ref{thm:derivative formula}).
Added smoothness should be an advantage in computational situations,
as was pointed out in \cite{Schreiber}.

It is not clear if the function in (\ref{eq:partion_three}) is continuous
whenever $\mu$ is finite with compact support.

The bound of the last partition function that we found in section
\ref{sec:Differentiating-the-Norm} used a different normalizing constant.
Recall we established that there is a $C$ so that \[
-\delta\leq\frac{d}{d\lambda}\ln\left(\left\Vert g_{e^{\lambda}}\ast\mu\right\Vert _{q}\right)\leq C.\]
In the Gaussian case, we had $C=0.$ Since \begin{eqnarray*}
 &  & \ln\left(\left(\int_{\mathbb{R}^{\delta}}\left(\int_{\mathbb{R}^{\delta}}g\left(\frac{\mathbf{x}-\mathbf{y}}{e^{\lambda}}\right)\, d\mu(\mathbf{y})\right)^{q}\frac{1}{e^{\delta\lambda}}\, dm(\mathbf{x})\right)^{\frac{1}{q-1}}\right)\\
 & = & \ln\left(e^{\delta\lambda}\left\Vert g_{e^{\lambda}}\ast\mu\right\Vert _{q}^{\frac{q}{q-1}}\right)\\
 & = & \delta\lambda+\frac{q}{q-1}\ln\left(\left\Vert g_{e^{\lambda}}\ast\mu\right\Vert _{q}\right)\end{eqnarray*}
we have\[
\frac{\delta}{1-q}\leq\frac{d}{d\lambda}\ln\left(\left(\int_{\mathbb{R}^{\delta}}\left(\int_{\mathbb{R}^{\delta}}g\left(\frac{\mathbf{x}-\mathbf{y}}{e^{\lambda}}\right)\, d\mu(\mathbf{y})\right)^{q}\frac{1}{e^{\delta\lambda}}\, dm(\mathbf{x})\right)^{\frac{1}{q-1}}\right)\leq C_{1}.\]
In the Gaussian case, we may take $C_{1}=\delta.$ 

In this section we prove that for $2\leq q<\infty,$ there is a constant
$C$ depending on $g$ and $q$ so that for any finite Borel measure
$\mu$ of compact support,\[
0\leq\frac{d}{d\lambda}\ln\left(\left(\int_{\mathbb{R}^{\delta}}\left(\int_{\mathbb{R}^{\delta}}g\left(\frac{\mathbf{x}-\mathbf{y}}{e^{\lambda}}\right)\, d\mu(\mathbf{y})\right)^{q-1}\, d\mu(\mathbf{x})\right)^{\frac{1}{q-1}}\right)\leq C,\]

We will need a lower bound on \[
\int_{\mathbb{R}^{\delta}}\epsilon^{\delta}\left(h_{\epsilon}\ast\mu\right)^{q-1}\, d\mu,\]
where $h$ is the negative of the radial derivative of $g.$ Since
$h(\mathbf{0})=0$ we need a small modification of the result in \cite{BarGerTch}.
We are restricting our attention to the case $1\leq q<\infty,$ which
allows us to avoid the technicalities encountered in \cite{BarGerTch}. 

Here we use the notation from \cite{Loring}, so $\mu^{(\epsilon)}$
is the sequence over $\mathbb{Z}^{\delta}$ given by\[
\mu_{\mathbf{n}}^{(\epsilon)}=\mu(\epsilon\mathbf{n}+\epsilon\mathbb{I}).\]

\begin{lem}
\label{lem:boundsOnCorrelConv}Assume that $g\geq0$ is rapidly decreasing
and that $1<q<\infty.$ There is a finite constant $C$ so that for
any finite Borel measure $\mu$ on $\mathbb{R}^{\delta},$\[
\frac{\epsilon^{\delta}\left\Vert g_{\epsilon}\ast\mu\right\Vert _{\mu,q-1}}{\left(S_{\mu}^{q}(\epsilon)\right)^{\frac{1}{q-1}}}\leq C\]
for all $\epsilon>0.$ If also $g(\mathbf{0})>0$ then there is a
$c>0$ so that\[
c\leq\frac{\epsilon^{\delta}\left\Vert g_{\epsilon}\ast\mu\right\Vert _{\mu,q-1}}{\left(S_{\mu}^{q}(\epsilon)\right)^{\frac{1}{q-1}}}.\]

\end{lem}
\begin{proof}
Define $\Gamma$ over $\mathbb{Z}^{\delta}$ by\[
\Gamma_{\mathbf{n}}=\sup\{ g(\mathbf{x})\mid\mathbf{x}\in\mathbf{n}+\mathbb{D}\},\]
where \[
\mathbb{D}=(-1,1)\times(-1,1)\times\cdots\times(-1,1).\]
\ifthenelse{\boolean{Details}}{If \[
\mathbf{x}\in\epsilon\mathbf{j}+\epsilon\mathbb{I}\]
and\[
\mathbf{y}\in\epsilon\mathbf{k}+\epsilon\mathbb{I}\]
then\[
\epsilon^{-1}(\mathbf{x}-\mathbf{y})\in(\mathbf{j}-\mathbf{k})+\mathbb{D}\]
and so\[
g(\epsilon^{-1}(\mathbf{x}-\mathbf{y}))\leq\Gamma_{\mathbf{\mathbf{j}-\mathbf{k}}}.\]

The assumptions of $g$ being rapidly decreasing and $\mu$ being
finite tells us\[
\left\Vert \Gamma\right\Vert _{1}<\infty.\]
}{}Repeating an argument from \cite{BarGerTch}, we find\begin{eqnarray*}
 &  & \left\Vert g_{\epsilon}\ast\mu\right\Vert _{\mu,q-1}^{q-1}\\
 & = & \epsilon^{-(q-1)\delta}\int\left(\int g(\epsilon^{-1}(\mathbf{x}-\mathbf{y}))\, d\mu(\mathbf{y})\right)^{q-1}\, d\mu(\mathbf{x})\\
 & = & \epsilon^{-(q-1)\delta}\sum_{\mathbf{j}\in\mathbb{Z}^{\delta}}\int_{\epsilon\mathbf{j}+\epsilon\mathbb{I}}\left(\sum_{\mathbf{k}\in\mathbb{Z}^{\delta}}\int_{\epsilon\mathbf{k}+\epsilon\mathbb{I}}g(\epsilon^{-1}(\mathbf{x}-\mathbf{y}))\, d\mu(\mathbf{y})\right)^{q-1}\, d\mu(\mathbf{x})\\
 & \leq & \epsilon^{-(q-1)\delta}\sum_{\mathbf{j}\in\mathbb{Z}^{\delta}}\left(\sum_{\mathbf{k}\in\mathbb{Z}^{\delta}}\Gamma_{\mathbf{\mathbf{j}-\mathbf{k}}}\mu_{\mathbf{k}}^{(\epsilon)}\right)^{q-1}\mu_{\mathbf{j}}^{(\epsilon)}.\end{eqnarray*}
H\"{o}lder's inequality and Young's convolution inequality now tell
us \begin{eqnarray*}
\left\Vert g_{\epsilon}\ast\mu\right\Vert _{\mu,q-1}^{q-1} & \leq & \epsilon^{-(q-1)\delta}\left(\left\Vert \Gamma\ast\mu^{(\epsilon)}\right\Vert _{q}\right)^{q-1}\left\Vert \mu^{(\epsilon)}\right\Vert _{q}\\
 & \leq & \epsilon^{-(q-1)\delta}\left\Vert \Gamma\right\Vert _{1}^{q-1}\left\Vert \mu^{(\epsilon)}\right\Vert _{q}^{q},\end{eqnarray*}
i.e.\[
\frac{\epsilon^{\delta}\left\Vert g_{\epsilon}\ast\mu\right\Vert _{\mu,q-1}}{\left(S_{\mu}^{p}(\epsilon)\right)^{\frac{1}{q-1}}}\leq\left\Vert \Gamma\right\Vert _{1}.\]

If $g$ is positive at the origin, then since it is continuous, we
can rescale $g$ using $\tilde{g}=g_{\eta}$ with the same properties
as $g,$ but with\[
\inf\{\tilde{g}(\mathbf{x})\mid\mathbf{x}\in\mathbb{D}\}>0\]
and\[
g_{\epsilon}=\tilde{g}_{\eta\epsilon}.\]
We can compare $\left\Vert g_{\epsilon}\ast\mu\right\Vert _{\mu,q-1}$
and $\left\Vert \tilde{g}_{\eta\epsilon}\ast\mu\right\Vert _{\mu,q-1}$
as follows:\begin{eqnarray*}
\frac{\epsilon^{\delta}\left\Vert g_{\epsilon}\ast\mu\right\Vert _{\mu,q-1}}{\left(S_{\mu}^{\delta}(\epsilon)\right)^{\frac{1}{q-1}}} & = & \frac{\epsilon^{\delta}\left\Vert \tilde{g}_{\eta\epsilon}\ast\mu\right\Vert _{\mu,q-1}}{\left(S_{\mu}^{\delta}(\epsilon)\right)^{\frac{1}{q-1}}}\\
 & = & \left(\eta^{-\delta}\right)\frac{\eta^{\delta}\epsilon^{\delta}\left\Vert \tilde{g}_{\eta\epsilon}\ast\mu\right\Vert _{\mu,q-1}}{\left(S_{\mu}^{q}(\eta\epsilon)\right)^{\frac{1}{q-1}}}\frac{\left(S_{\mu}^{q}(\eta\epsilon)\right)^{\frac{1}{q-1}}}{\left(S_{\mu}^{q}(\epsilon)\right)^{\frac{1}{q-1}}}\end{eqnarray*}
By \cite[Theorem 3.4]{Loring}, there are constants $A$ and $B$
so that \[
e^{-A-B|\ln(\eta)|}\leq\frac{S_{\mu}^{q}(\eta\epsilon)}{S_{\mu}^{q}(\epsilon)}\leq e^{A+B|\ln(\eta)|}\]
for all $\epsilon.$ Therefore\[
\frac{\epsilon^{\delta}\left\Vert g_{\epsilon}\ast\mu\right\Vert _{\mu,q-1}}{\left(S_{\mu}^{q}(\epsilon)\right)^{\frac{1}{q-1}}}\leq\left(\eta^{-\delta}e^{\frac{A+B|\ln(\eta)|}{q-1}}\right)\frac{(\eta\epsilon)^{\delta}\left\Vert \tilde{g}_{\eta\epsilon}\ast\mu\right\Vert _{\mu,q-1}}{\left(S_{\mu}^{q}(\eta\epsilon)\right)^{\frac{1}{q-1}}},\]
and it suffices to prove the result in the case where\[
\inf\{ g(\mathbf{x})\mid\mathbf{x}\in\mathbb{D}\}>0.\]

Let\[
\gamma_{\mathbf{n}}=\inf\{ g(\mathbf{x})\mid\mathbf{x}\in\mathbf{n}+\mathbb{D}\}.\]
As above, we find\[
\left\Vert g_{\epsilon}\ast\mu\right\Vert _{\mu,q-1}^{q-1}\geq\epsilon^{-(q-1)\delta}\sum_{\mathbf{j}\in\mathbb{Z}^{\delta}}\left(\sum_{\mathbf{k}\in\mathbb{Z}^{\delta}}\gamma_{\mathbf{\mathbf{j}-\mathbf{k}}}\mu_{\mathbf{k}}^{(\epsilon)}\right)^{q-1}\mu_{\mathbf{j}}^{(\epsilon)}\]
and so\begin{eqnarray*}
\left\Vert g_{\epsilon}\ast\mu\right\Vert _{\mu,q-1}^{q-1} & \geq & \epsilon^{-(q-1)\delta}\sum_{\mathbf{j}\in\mathbb{Z}^{\delta}}\left(\gamma_{\mathbf{0}}\mu_{\mathbf{j}}^{(\epsilon)}\right)^{q-1}\mu_{\mathbf{j}}^{(\epsilon)}\\
 & = & \gamma_{\mathbf{0}}\epsilon^{-(q-1)\delta}\sum_{\mathbf{j}\in\mathbb{Z}^{\delta}}\left(\mu_{\mathbf{j}}^{(\epsilon)}\right)^{q}\end{eqnarray*}
and \[
\epsilon^{\delta}\left\Vert g_{\epsilon}\ast\mu\right\Vert _{\mu,q-1}\geq\gamma_{\mathbf{0}}^{\frac{1}{q-1}}\left(S_{\mu}^{q}(\epsilon)\right)^{\frac{1}{q-1}}.\]

\end{proof}
\begin{thm}
\label{thm:derivative formula 2}Suppose $2\leq q<\infty.$ Suppose
$g:\mathbb{R}^{\delta}\rightarrow\mathbb{R}$ is rapidly decreasing
and let $h:\mathbb{R}^{\delta}\rightarrow\mathbb{R}$ denote the negative
of the radial derivative of $g.$ Suppose $g\geq0$ and $h\geq0.$
If $\mu$ is a finite Borel measure on $\mathbb{R}^{\delta}$ then
\[
\frac{d}{d\epsilon}\left\Vert g_{\epsilon}\ast\mu\right\Vert _{\mu,q-1}=\frac{\int_{\mathbb{R}^{\delta}}\left(g_{\epsilon}\ast\mu\right)^{q-2}h_{\epsilon}\ast\mu\, d\mu}{\epsilon\left\Vert g_{\epsilon}\ast\mu\right\Vert _{\mu,q-1}^{q-2}}-\frac{\delta\left\Vert g_{\epsilon}\ast\mu\right\Vert _{\mu,q-1}}{\epsilon}\]
and\[
\frac{d}{d\lambda}\ln\left(\left\Vert g_{e^{\lambda}}\ast\mu\right\Vert _{\mu,q-1}\right)=\frac{\int_{\mathbb{R}^{\delta}}\left(g_{e^{\lambda}}\ast\mu\right)^{q-2}h_{e^{\lambda}}\ast\mu\, d\mu}{\int_{\mathbb{R}^{\delta}}\left(g_{e^{\lambda}}\ast\mu\right)^{q-1}\, d\mu}-\delta.\]

\end{thm}
\begin{proof}
\ifthenelse{\boolean{Details}}{Assume $0<a\leq\epsilon\leq b.$

By Lemma \ref{lem:inner derivative formula},

\begin{eqnarray*}
\left|\frac{\partial}{\partial\epsilon}\left[(g_{\epsilon}\ast\mu)(\mathbf{x})\right]\right| & = & \left|\epsilon^{-1}\left((h_{\epsilon}\ast\mu)(\mathbf{x})-\delta(g_{\epsilon}\ast\mu)(\mathbf{x})\right)\right|\\
 & \leq & \epsilon^{-1}\left(\left\Vert \mu\right\Vert _{1}\left\Vert h_{\epsilon}\right\Vert _{\infty}+\delta\left\Vert \mu\right\Vert _{1}\left\Vert g_{\epsilon}\right\Vert _{\infty}\right)\\
 & = & \epsilon^{-\delta-1}\left\Vert \mu\right\Vert _{1}\left(\left\Vert h\right\Vert _{\infty}+\delta\left\Vert g\right\Vert _{\infty}\right)\\
 & \leq & a^{-\delta-1}\left\Vert \mu\right\Vert _{1}\left(\left\Vert h\right\Vert _{\infty}+\delta\left\Vert g\right\Vert _{\infty}\right).\end{eqnarray*}
Also\begin{eqnarray*}
\left((g_{\epsilon}\ast\mu)(\mathbf{x})\right)^{q-2} & \leq & a^{-(q-2)(\delta+1)}\left\Vert \mu\right\Vert _{1}^{q-2}\left\Vert g\right\Vert _{\infty}^{q-2}\end{eqnarray*}
and so}{For $\epsilon$ restricted to some interval $[a,b],$ it
follows from Lemma \ref{lem:inner derivative formula} that}\[
\left|\frac{\partial}{\partial\epsilon}\left((g_{\epsilon}\ast\mu)(\mathbf{x})\right)^{q-1}\right|\leq(q-1)a^{-(q-1)(\delta+1)}\left\Vert \mu\right\Vert _{1}^{q-1}\left\Vert g\right\Vert _{\infty}^{q-2}\left(\left\Vert h\right\Vert _{\infty}+\delta\left\Vert g\right\Vert _{\infty}\right).\]
Dominated convergence yields\[
\frac{d}{d\epsilon}\left(\int_{\mathbb{R}^{\delta}}\left(g_{\epsilon}\ast\mu\right)^{q-1}\, d\mu\right)=\frac{q-1}{\epsilon}\left(\int_{\mathbb{R}^{\delta}}\left(g_{\epsilon}\ast\mu\right)^{q-2}h_{\epsilon}\ast\mu\, d\mu-\delta\left\Vert g_{\epsilon}\ast\mu\right\Vert _{\mu,q-1}^{q-1}\right)\]
 and so\begin{eqnarray*}
 &  & \frac{d}{d\epsilon}\left\Vert g_{\epsilon}\ast\mu\right\Vert _{\mu,q-1}\\
 & = & \frac{1}{\epsilon}\left\Vert g_{\epsilon}\ast\mu\right\Vert _{\mu,q-1}^{2-q}\left(\int_{\mathbb{R}^{d}}\left(g_{\epsilon}\ast\mu\right)^{q-2}h_{\epsilon}\ast\mu\, d\mu-\delta\left\Vert g_{\epsilon}\ast\mu\right\Vert _{\mu,q-1}^{q-1}\right)\\
 & = & \frac{\int_{\mathbb{R}^{\delta}}\left(g_{\epsilon}\ast\mu\right)^{q-2}h_{\epsilon}\ast\mu\, d\mu}{\epsilon\left\Vert g_{\epsilon}\ast\mu\right\Vert _{\mu,q-1}^{q-2}}-\frac{\delta\left\Vert g_{\epsilon}\ast\mu\right\Vert _{\mu,q-1}}{\epsilon}.\end{eqnarray*}
We use\[
\frac{d}{d\lambda}\ln\left(\left\Vert g_{e^{\lambda}}\ast\mu\right\Vert _{\mu,q-1}\right)=\frac{\left(\left.\frac{d}{d\epsilon}\right|_{\epsilon=e^{\lambda}}\left\Vert g_{\epsilon}\ast\mu\right\Vert _{\mu,q-1}\right)e^{\lambda}}{\left\Vert g_{e^{\lambda}}\ast\mu\right\Vert _{\mu,q-1}}\]
and find\[
\frac{d}{d\lambda}\ln\left(\left\Vert g_{e^{\lambda}}\ast\mu\right\Vert _{\mu,q-1}\right)=\frac{\int_{\mathbb{R}^{\delta}}\left(g_{e^{\lambda}}\ast\mu\right)^{q-2}h_{e^{\lambda}}\ast\mu\, d\mu}{\left\Vert g_{e^{\lambda}}\ast\mu\right\Vert _{\mu,q-1}^{q-1}}-\delta.\]

\end{proof}
\begin{thm}
\label{thm:boundCorrIntDeriv}Suppose $2\leq q<\infty.$ Suppose $g:\mathbb{R}^{\delta}\rightarrow\mathbb{R}$
is rapidly decreasing and let $h:\mathbb{R}^{\delta}\rightarrow\mathbb{R}$
denote the negative of the radial derivative of $g.$ Suppose $g\geq0$
and $h\geq0.$ If $\mu$ is a finite Borel measure of compact support
on $\mathbb{R}^{d}$ then\[
-\delta\epsilon^{-1}\left\Vert g_{\epsilon}\ast\mu\right\Vert _{\mu,q-1}\leq\frac{d}{d\epsilon}\left(\left\Vert g_{\epsilon}\ast\mu\right\Vert _{\mu,q-1}\right)\leq\epsilon^{-1}\left\Vert h_{\epsilon}\ast\mu\right\Vert _{\mu,q-1}-\delta\epsilon^{-1}\left\Vert g_{\epsilon}\ast\mu\right\Vert _{\mu,q-1}\]
and\[
-\delta\leq\frac{d}{d\lambda}\ln\left(\left\Vert g_{e^{\lambda}}\ast\mu\right\Vert _{\mu,q-1}\right)\leq\frac{\left\Vert h_{e^{\lambda}}\ast\mu\right\Vert _{\mu,q-1}}{\left\Vert g_{e^{\lambda}}\ast\mu\right\Vert _{\mu,q-1}}-\delta.\]

\end{thm}
\begin{proof}
If $2<q<\infty,$ we can apply H\"{o}lder's inequality and we find\[
\int_{\mathbb{R}^{\delta}}\left(g_{\epsilon}\ast\mu\right)^{q-2}\left(h_{\epsilon}\ast\mu\right)\, d\mu\leq\left(\int_{\mathbb{R}^{\delta}}\left(g_{\epsilon}\ast\mu\right)^{q-1}\, d\mu\right)^{\frac{q-2}{q-1}}\left(\int_{\mathbb{R}^{\delta}}\left(h_{\epsilon}\ast\mu\right)^{q-1}\, d\mu\right)^{\frac{1}{q-1}}.\]
This is trivially true as well when $q=2.$ We can rewrite this as\[
\int_{\mathbb{R}^{\delta}}\left(g_{\epsilon}\ast\mu\right)^{q-2}\left(h_{\epsilon}\ast\mu\right)\, d\mu\leq\left\Vert g_{\epsilon}\ast\mu\right\Vert _{\mu,q-1}^{q-2}\left\Vert h_{\epsilon}\ast\mu\right\Vert _{\mu,q-1}\]
and the inequalities follow from the last result and the fact that
$g_{\epsilon}\ast\mu$ and $h_{\epsilon}\ast\mu$ are nonnegative.
\end{proof}
\begin{cor}
Suppose $2\leq q<\infty.$ Suppose $g:\mathbb{R}^{\delta}\rightarrow\mathbb{R}$
is nonnegative, nontrivial, rapidly decreasing and is radially nonincreasing.
There is a finite constant $C$ so that if $\mu$ is a finite Borel
measure of compact support on $\mathbb{R}^{\delta}$ then \[
0\leq\frac{d}{d\lambda}\ln\left(\left(\int_{\mathbb{R}^{\delta}}\left(\int_{\mathbb{R}^{\delta}}g\left(\frac{\mathbf{x}-\mathbf{y}}{e^{\lambda}}\right)\, d\mu(\mathbf{y})\right)^{q-1}\, d\mu(\mathbf{x})\right)^{\frac{1}{q-1}}\right)\leq C.\]
 
\end{cor}
\begin{proof}
By Lemma \ref{lem:boundsOnCorrelConv}, we have an upper bound on
$\left\Vert h_{e^{\lambda}}\ast\mu\right\Vert _{\mu,q-1}$ and a lower
bound on $\left\Vert g_{e^{\lambda}}\ast\mu\right\Vert _{\mu,q-1}$
that depends only on $q$ and $g.$ Therefore Theorem \ref{thm:continuityCorrelationInt}
gives us a $C_{1}$ so that\[
-\delta\leq\frac{d}{d\lambda}\ln\left(\left\Vert g_{e^{\lambda}}\ast\mu\right\Vert _{\mu,q-1}\right)\leq C_{1}\]
for all $\mu$ and all $\lambda.$ As to the partition function,\[
\ln\left(\left(\int_{\mathbb{R}^{\delta}}\left(\int_{\mathbb{R}^{\delta}}g\left(\frac{\mathbf{x}-\mathbf{y}}{e^{\lambda}}\right)\, d\mu(\mathbf{y})\right)^{q-1}\, d\mu(\mathbf{x})\right)^{\frac{1}{q-1}}\right)=\delta\lambda+\ln\left(\left\Vert g_{e^{\lambda}}\ast\mu\right\Vert _{\mu,q-1}\right)\]
and so\[
0\leq\frac{d}{d\lambda}\ln\left(\left(\int_{\mathbb{R}^{\delta}}\left(\int_{\mathbb{R}^{\delta}}g\left(\frac{\mathbf{x}-\mathbf{y}}{e^{\lambda}}\right)\, d\mu(\mathbf{y})\right)^{q-1}\, d\mu(\mathbf{x})\right)^{\frac{1}{q-1}}\right)\leq C_{1}+\delta.\]

\end{proof}
\begin{thm}
\label{thm:continuityCorrelationInt} Assume that $g\geq0$ is rapidly
decreasing and that $1<q<\infty.$ For any finite Borel measure $\mu$
on $\mathbb{R}^{\delta}$ with compact support, \[
\int_{\mathbb{R}^{\delta}}\left(g_{\epsilon}\ast\mu\right)^{q-1}\, d\mu(\mathbf{x})\]
varies continuously in $\epsilon.$
\end{thm}
\begin{proof}
Assume $0<a\leq\epsilon\leq b.$ If $G$ is a bound on $g,$ then
$a^{-\delta}G$ is a bound on $g_{\epsilon}$ and so\[
\left(\left(g_{\epsilon}\ast\mu\right)(\mathbf{y})\right)\leq a^{-\delta(q-1)}G^{q-1}\|\mu\|_{1}^{q-1}.\]
Since $\mu$ is a finite measure, we can apply the dominated convergence
theorem and have\[
\lim_{\epsilon\rightarrow\eta}\int_{\mathbb{R}^{\delta}}\left(g_{\epsilon}\ast\mu\right)^{q-1}\, d\mu(\mathbf{x})=\int_{\mathbb{R}^{\delta}}\left(g_{\eta}\ast\mu\right)^{q-1}\, d\mu(\mathbf{x}).\]
 
\end{proof}

\section{Gaussian kernel R\'{e}nyi entropy sums\label{sec:Gaussian-Kernel-Sums}}

There is a smooth version of the partition function\[
\sum_{\mathbf{j}\in\mathbb{Z}^{\delta}}\mu(\epsilon\mathbf{j}+\epsilon\mathbb{I})^{q}\]
that eliminates the sharp cut-off at the boundary of the cells in
the grid,\[
\sum_{\mathbf{j}\in\mathbb{Z}^{\delta}}\left(\int_{\mathbb{R}^{\delta}}g\left(\mathbf{j}-\frac{\mathbf{y}}{\epsilon}\right)\, d\mu(\mathbf{y})\right)^{q}.\]
Although we have not determined if this creates a partition function
that is differentiable, it does at least give continuity.

What we have in mind for $g$ is either a Gaussian, or a smooth function
between $0$ and $1$ that equals the characteristic function for
$\mathbb{I}$ except close to the boundary of $\mathbb{I}.$

Recall from \cite{Loring} that we say a finite Borel measure $\mu$
on $\mathbb{R}^{d}$ is $q$\emph{-finite} if $S_{\mu}^{q}(1)<\infty.$
This is automatic if $1<q<\infty.$ 

First we show that this modified R\'{e}nyi entropy sum still leads
to $D_{q}^{\pm}(\mu).$

\begin{thm}
Assume that $g\geq0$ is rapidly decreasing, with $g(\mathbf{0})>0,$
and that $0<q<\infty,$ $q\neq1.$ There is a constant $C$ so that,
for any Borel measure $\mu$ on $\mathbb{R}^{\delta},$ \[
C^{-1}\leq\frac{\sum_{\mathbf{j}\in\mathbb{Z}^{\delta}}\left(\int_{\mathbb{R}^{\delta}}g\left(\mathbf{j}-\frac{\mathbf{y}}{\epsilon}\right)\, d\mu(\mathbf{y})\right)^{q}}{\sum_{\mathbf{j}\in\mathbb{Z}^{\delta}}\mu(\epsilon\mathbf{j}+\epsilon\mathbb{I})^{q}}\leq C\]
for all $\epsilon>0.$
\end{thm}
\begin{proof}
The proof is almost identical to that of \cite[Lemma 2.3]{Loring},
and we again use the notation used there.

Notice that\[
\sum_{\mathbf{j}\in\mathbb{Z}^{\delta}}\left(\int_{\mathbb{R}^{\delta}}g\left(\mathbf{j}-\frac{\mathbf{y}}{\epsilon}\right)\, d\mu(\mathbf{y})\right)^{q}=\sum_{\mathbf{j}\in\mathbb{Z}^{\delta}}\left(\sum_{\mathbf{k}\in\mathbb{Z}^{\delta}}\int_{\epsilon\mathbf{k}+\epsilon\mathbb{I}}g\left(\mathbf{j}-\frac{\mathbf{y}}{\epsilon}\right)\, d\mu(\mathbf{y})\right)^{q}.\]
If \[
\mathbf{y}\in\epsilon\mathbf{k}+\epsilon\mathbb{I}\]
then\[
\mathbf{j}-\frac{\mathbf{y}}{\epsilon}\in(\mathbf{j}-\mathbf{k})-\mathbb{I}\subseteq(\mathbf{j}-\mathbf{k})-\mathbb{D}.\]
Therefore\[
\sum_{\mathbf{j}\in\mathbb{Z}^{\delta}}\left(\int_{\mathbb{R}^{\delta}}g\left(\mathbf{j}-\frac{\mathbf{y}}{\epsilon}\right)\, d\mu(\mathbf{y})\right)^{q}\leq\sum_{\mathbf{j}\in\mathbb{Z}^{d}}\left(\sum_{\mathbf{k}\in\mathbb{Z}^{d}}\Gamma_{\mathbf{j}-\mathbf{k}}\mu(\epsilon\mathbf{k}+\epsilon\mathbb{I})\right)^{q}\]
and\[
\sum_{\mathbf{j}\in\mathbb{Z}^{\delta}}\left(\int_{\mathbb{R}^{\delta}}g\left(\mathbf{j}-\frac{\mathbf{y}}{\epsilon}\right)\, d\mu(\mathbf{y})\right)^{q}\geq\sum_{\mathbf{j}\in\mathbb{Z}^{d}}\left(\sum_{\mathbf{k}\in\mathbb{Z}^{d}}\gamma_{\mathbf{j}-\mathbf{k}}\mu(\epsilon\mathbf{k}+\epsilon\mathbb{I})\right)^{q}.\]
Therefore\[
\left\Vert \gamma\ast\mu^{(\epsilon)}\right\Vert _{q}^{q}\leq\sum_{\mathbf{j}\in\mathbb{Z}^{\delta}}\left(\int_{\mathbb{R}^{\delta}}g\left(\mathbf{j}-\frac{\mathbf{y}}{\epsilon}\right)\, d\mu(\mathbf{y})\right)^{q}\leq\left\Vert \Gamma\ast\mu^{(\epsilon)}\right\Vert _{q}^{q}\]
and the rest of the proof follows that of \cite[Lemma 2.3]{Loring}.
\end{proof}
\begin{thm}
\label{thm:continuityGaussKerSums}Assume that $g\geq0$ is rapidly
decreasing and that $0<q<\infty,$ $q\neq1.$ For any finite Borel
measure $\mu$ on $\mathbb{R}^{\delta}$ with compact support, the
sum \[
\sum_{\mathbf{j}\in\mathbb{Z}^{\delta}}\left(\int_{\mathbb{R}^{\delta}}g\left(\mathbf{j}-\frac{\mathbf{y}}{\epsilon}\right)\, d\mu(\mathbf{y})\right)^{q}\]
varies continuously in $\epsilon.$
\end{thm}
\begin{proof}
Since\[
\sum_{\mathbf{j}\in\mathbb{Z}^{\delta}}\left(\int_{\mathbb{R}^{\delta}}g\left(\mathbf{j}-\frac{\mathbf{y}}{\epsilon}\right)\, d\mu(\mathbf{y})\right)^{q}=\epsilon^{\delta q}\sum_{\mathbf{j}\in\mathbb{Z}^{\delta}}\left(g_{\epsilon}\ast\mu(\epsilon\mathbf{j})\right)^{q}\]
it will suffice to prove the continuity of\[
\sum_{\mathbf{j}\in\mathbb{Z}^{\delta}}\left(g_{\epsilon}\ast\mu(\epsilon\mathbf{j})\right)^{q}.\]
Again we restrict $\epsilon$ to some the interval $[a,b].$ 

We know that $g$ is bounded by some $G<\infty.$ We have the bound\[
\epsilon^{-\delta}g\left(\mathbf{j}-\frac{\mathbf{y}}{\epsilon}\right)\leq a^{-\delta}G,\]
and since $\mu$ is a finite measure, we can apply dominated convergence
and conclude that

\[
g_{\epsilon}\ast\mu(\epsilon\mathbf{j})=\int_{\mathbb{R}^{\delta}}\epsilon^{-\delta}g\left(\mathbf{j}-\frac{\mathbf{y}}{\epsilon}\right)\, d\mu(\mathbf{y})\]
 varies continuously in $\epsilon.$

We saw in the proof of Theorem \ref{thm:derivative formula} that
there is a bound\[
(g_{\epsilon}\ast\mu)(\mathbf{x})\leq C_{2}(C_{3}\wedge|\mathbf{x}|^{-k}),\]
where $k$ is taken to be an integer larger than $\frac{\delta+1}{q}.$
Therefore\[
\left((g_{\epsilon}\ast\mu)(\epsilon\mathbf{j})\right)^{q}\leq a^{-qk}C_{2}(b^{qk}C_{3}^{q}\wedge|\mathbf{j}|^{-qk}).\]
We can apply dominated convergence, this time with the measure being
counting measure on $\mathbb{Z}^{\delta},$ and conclude that\[
\sum_{\mathbf{j}\in\mathbb{Z}^{\delta}}\left(g_{\epsilon}\ast\mu(\epsilon\mathbf{j})\right)^{q}\]
varies continuously in $\epsilon$.
\end{proof}

\end{document}